\documentclass[a4paper]{article}    
\usepackage{amsthm,amssymb,amsmath,enumerate,graphicx}

\theoremstyle{plain}
\newtheorem{definition}{Definition}[section]

\newtheorem{theorem}[definition]{Theorem}

\theoremstyle{definition}
\newtheorem{example}[definition]{Example}

\newcommand{\fin}{\rm fin}

\newcommand{\CCC}{\mathcal{C}}
\newcommand{\CCCfin}{\mathcal{C}_{\fin}}

\newcommand{\EEE}{\mathcal{E}}
\newcommand{\FF}{\mathbb{F}}

\newcommand{\NN}{\mathbb{N}}

\newcommand{\ZZ}{\mathbb{Z}}

\newcommand{\sm}{\setminus}
\renewcommand{\ker}{{\rm Ker\>}}
\newcommand{\im}{{\rm Im\>}}
\newcommand{\kreis}[1]{\mathaccent"7017\relax #1}

\newcommand{\restr}{\!\restriction\!}

\newcommand{\es}{\emptyset}

\newcommand{\fwd}[3]{%
  \overset{\smash{\raisebox{-#2pt}{\hspace*{#3 pt}\tiny$\rightarrow$}}}{#1}}
\newcommand{\bkwd}[3]{%
  \overset{\smash{\raisebox{-#2pt}{\hspace*{#3 pt}\tiny$\leftarrow$}}}{#1}}
  
\newcommand{\ve}[1][\hspace*{1pt}]{\fwd{e}{1.75}{1}_{\hspace*{-1pt}#1}}
\newcommand{\ev}[1][\hspace*{1pt}]{\bkwd{e}{1.75}{0}_{\hspace*{-1pt}#1}}

\newcommand{\vCCC}[1][\hspace*{1pt}]{\fwd{\CCC}{1.75}{1}_{\hspace*{-1pt}#1}}

\newcommand{\COMMENT}[1]{}

\newcommand{\assign}{
  \mathrel{\mathop{:}}=
}

\title{\hbox{\llap{Lo}cally finite graphs with ends: a topological approach.}
  III.~Fundamental group and homology}

\author{Reinhard Diestel and Philipp Spr\"ussel}
\date{}

\begin{document}      

\maketitle

\begin{abstract}\noindent
  This paper is the last part of a comprehensive survey of a newly
  emerging field: a topological approach to the study of locally
  finite graphs that crucially incorporates their ends. Topological
  arcs and circles, which may pass through ends, assume the role
  played in finite graphs by paths and cycles. The first two parts of
  the survey together provide a suitable entry point to this field
  for new readers; they are available in combined form from the
  ArXiv~\cite{RDsBanffSurvey}.
  
  The topological approach indicated above has made it possible to
  extend to locally finite graphs many classical theorems of finite
  graph theory that do not extend verbatim. While the first
  part~\cite{TopSurveyI} of this survey introduces the theory as such
  and the second part~\cite{TopSurveyII} is devoted to those
  applications, this third part looks at the theory from an algebraic-topological point
  of view.

The results surveyed here include both a combinatorial description of the fundamental group of a locally finite graph with ends and the homology aspects of this space.
\end{abstract}

\section{Introduction}\label{sec:intro}

The survey~\cite{RDsBanffSurvey} describes a topological framework in
which many well-known theorems about finite graphs that appear to
fail for infinite graphs do have a natural infinite analogue. It has
been realised in recent years that many such theorems, especially
about paths and cycles, work in a slightly richer setting: not in
the (locally finite) graph $G$ itself, but in its compactification
$|G|$ obtained by adding its \emph{ends}.%
  \footnote{For a formal definition of $|G|$ see~\cite{RDsBanffSurvey}.}
In this setting, the traditional cycle space of a graph is replaced by
its \emph{topological cycle space}. The topological cycle space
$\CCC=\CCC(G)$ of a locally finite graph $G$ is based on (the edge
sets of) \emph{topological circles} in~$|G|$, homeomorphic images of
the unit circle $S^1$, allowing infinite sums as long as they are
\emph{thin}, that is, every edge appears in only finitely many
summands. Since the topological cycle space $\CCC(G)$ was introduced~\cite{CyclesI, CyclesII}, it has proved surprisingly successful; see~\cite{RDsBanffSurvey, TopSurveyII} for numerous applications.
   \COMMENT{}

Given the success of $\CCC$ for graphs, it seems desirable to recast
its definition in homological terms that make no reference to the
one-dimensional character of~$|G|$ (e.g., to circles), to obtain a
homology theory for similar but more general spaces (such as
non-compact CW complexes of any dimension) that implements the ideas
and advantages of $\CCC$ more generally. This approach has been
pursued in~\cite{FundGp,Hom1,Hom2}. In this paper we present its
main ideas, results and examples. For simplicity, all our coefficients will be taken from~$\FF_2$.

For such an extendable translation of our combinatorial definition of
$\CCC$ into algebraic terms, simplicial homology is easily seen not to
be the right approach: while $|G|$ is not a simplicial complex, the
simplicial homology of $G$ itself (without ends) yields the classical
cycle space $\CCCfin$. One way of extending simplicial homology to
more general spaces is \v{C}ech homology; and indeed we will show that
its first group applied to~$|G|$ is isomorphic to $\CCC$. But there
the usefulness of \v{C}ech homology for graphs ends: since its groups
are constructed as limits rather than directly from chains and cycles,
they do not interact with the combinatorial structure of $G$ in the
way we expect and know it from~$\CCC$.

The next candidate for the desired description of $\CCC$ in terms of
homology is singular homology. Indeed, $\CCC$ is built from circles in~$|G|$, and circles are singular $1$-cycles that generate the first
singular homology group $H_1(|G|)$ of $|G|$, so both groups are built
from similar elements. On the face of it, it is not clear whether
$\CCC$ might in fact be isomorphic, even canonically, to $H_1(|G|)$.
However, it will turn out that it is not: in~\cite{Hom1} we prove that
$\CCC$ is always a natural quotient of $H_1(|G|)$, and this quotient
is proper unless $G$ is essentially finite. This may seem surprising,
since $\CCC$ is defined via (thin) infinite sums while all sums in the
definition of $H_1(|G|)$ are finite, which suggests that $\CCC$ might
be larger than $H_1(|G|)$.

Our approach for the comparison of $\CCC$ and $H_1(|G|)$ will be to
define a homomorphism from $Z_1(|G|)$ to the edge space $\EEE$ that
counts how often the edges of $G$ are traversed by the simplices
of a $1$-cycle~$z$, and maps $z$ to the set of those edges that
are traversed an odd number of times. It will turn out that this
homomorphism vanishes on boundaries and that its image is precisely
$\CCC$. Hence it defines an epimorphism $f\colon H_1(|G|)\to\CCC(G)$.
However, we will show that $f$ is not normally injective. Indeed, there will be loops that traverse every edge evenly often (even equally
often in either direction), but which can be shown not with some effort to be null-homologous. Thus,
$\CCC$ is a genuinely new object, also from a topological point of
view.

For our proof that those loops are not
null-homologous we shall need a better understanding of the
fundamental group of~$|G|$. This will enable us to define an invariant
on $1$-chains in $|G|$ that can distinguish certain $1$-cycles from
boundaries of singular $2$-chains, hence completing the proof that $f$
need not be injective. The fundamental group of a finite graph $G$ is
easy to describe: it is the free group on the (oriented) \emph{chords}
of a spanning tree of $G$, the edges of $G$ that are not edges of the
spanning tree. For the Freudenthal compactification of infinite
graphs, the situation is different, since a loop in $|G|$ can traverse
infinitely many chords while the elements of a free group are always
finite sums of its generators.

One of the main aims of this project, therefore, became to develop a
combinatorial description of the fundamental group of the space~$|G|$
for an arbitrary connected locally finite graph~$G$. In~\cite{FundGp}
we describe~$\pi_1(|G|)$, as for finite~$G$, in terms of reduced words
in the oriented chords of a spanning tree. However, when $G$ is
infinite this does not work with arbitrary spanning trees but only with
\emph{topological spanning trees}. Moreover, we will have to allow
infinite words of any countable order type, and likewise allow the reduction sequences
cancelling adjacent inverse letters to have arbitrary countable order type.
However, these reductions can also be described in
terms of word reductions in the free groups $F_I$ on all the finite
subsets $I$ of chords, which enables us to embed the group $F_\infty$
of infinite reduced words as a subgroup in the inverse limit of those~$F_I$, and
handle it in this form. On the other hand, mapping a loop in $|G|$ to
the sequence of chords it traverses, and then reducing that sequence
(or word), turns out to be well defined on homotopy classes and hence
defines an embedding of $\pi_1(|G|)$ as a subgroup in~$F_\infty$.%
   \COMMENT{}

Having proved that $\CCC$ is usually a proper quotient of $H_1(|G|)$,
the last aim of this project then was to define a variant of singular
homology that works in more general spaces, and which for graphs captures precisely~$\CCC$. First steps in this direction were taken in~\cite{Hom1}; it was completed
in~\cite{Hom2}. Our hope with this translation was to stimulate further work
in two directions. One is that its new topological guise should make
the cycle space accessible to topological methods that might generate
some windfall for the study of graphs. And conversely, that as the
approach that gave rise to~$\CCC$ is made accessible for more general
spaces---in particular, for CW complexes of higher dimensions---its proven usefulness for graphs might find some more general topological analogues.

The key to the definition of~$\CCC$, and to its success, is that it treats ends differently from other points. To preserve this feature, our new homology theory is constructed for locally
compact Hausdorff spaces $X$ with a fixed Hausdorff compactification~$\hat X$, in which the compacification points play the role of ends.

\section{\v{C}ech homology}

The \v{C}ech homology of a space is an alternative to singular
homology for spaces that are not simplicial complexes. For a general
space $X$, the $n$th \v{C}ech homology group $\check H_n(X)$ is the
inverse limit of the homology groups of simplicial complexes induced
by open covers of $X$.%
  \footnote{See~\cite{Hom1} for a formal definition.}
In the case of $X=|G|$, one can compute the groups $\check H_n(X)$
more directly. To do so, fix a normal spanning tree $T$ of $G$, with root $r$
say, and denote the subtree of $T$ induced by the first $i$ levels by
$T_i$. Let $G_i$ be the finite graph obtained from $G$ by contracting
each component of $G-T_i$; then $\check H_n(X)$ is the inverse limit
of the family $\big(H_n(G_i),\le\big)_{i\in\NN}$. Since $\CCC(G)$ is
the inverse limit of the groups $H_1(G_i)$, we have

\begin{theorem}[\cite{Hom1}]\label{thm:CechequalsC}
  For a locally finite graph $G$ we have a canonical isomorphism
  $\check H_1(|G|) \simeq \CCC(G)$.
\end{theorem}

Theorem~\ref{thm:CechequalsC} shows that one can describe the
topological cycle space in terms of the \v{C}ech homology. However,
although $\check H_1(|G|)$ is isomorphic to $\CCC(G)$ as a group, it
does not sufficiently reflect the combinatorial properties of
$\CCC(G)$, its interaction with the combinatorial structure of $G$.
To make this precise, note that a number of classical results about
the cycle space say which circuits generate it---as do the
non-separating chordless circuits in a $3$-connected graph, say. In
the \v{C}ech homology, however, it is not possible to decide whether a
given homology class in $\check H_1(|G|)$ corresponds to a circuit. Indeed, the obvious relation between $\check H_1(|G|)$ and the
combinatorial structure of $G$ is that every homology class $c\in
\check H_1(|G|)$ corresponds to a family $(c_n)$ of homology classes
in the groups $H_1(G_n)$. One might think that the class $c$ should
correspond to a circuit in $|G|$ if and only if every $c_n$ with sufficiently
large $n$ corresponds to a circuit in $G_n$. But this is not the case:
the limit of a sequence of cycle space elements in the $G_n$ can be a
circuit even if the elements of the sequence are not circuits in
the~$G_n$.\COMMENT{}

\begin{figure}[htbp]
  \centering
  \includegraphics[width=.75\linewidth]{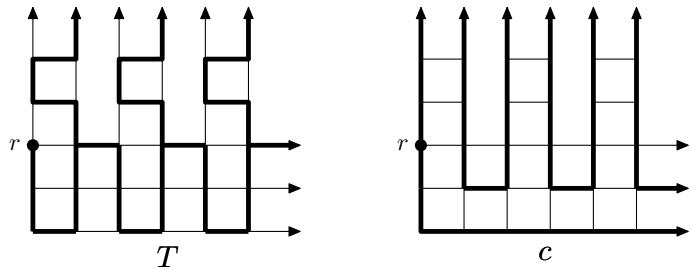}
  \caption{The graph $G$ (drawn twice) with a normal spanning tree $T$
  and a circuit $c$.}
  \label{fig:noncircuits:graph}
\end{figure}

Let $G$ be the graph shown in Figure~\ref{fig:noncircuits:graph}. $G$
consists of a `wide ladder' with three `poles' $x_1^1,x_2^1,\dotsc$,
$x_1^2,x_2^2,\dotsc$, and $x_1^3,x_2^3,\dotsc$, and has attached
infinitely many (oridinary) ladders by identifying the first rung of
the $n$th ladder $L_n$ with the edge $x_{2n-1}^1x_{2n}^1$. It is not
hard to prove that $T$ from Figure~\ref{fig:noncircuits:graph} is a
normal spanning tree of $G$ with root $r=x^1_1$. 

\begin{figure}[htbp]
  \centering
  \includegraphics[width=.75\linewidth]{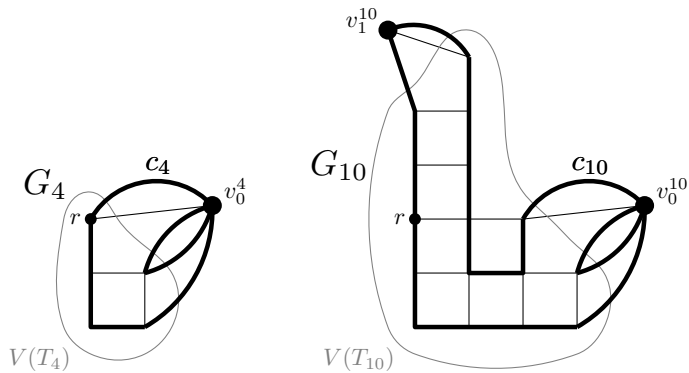}
  \caption{The edge sets $c_4$ in $G_4$ and $c_{10}$ in $G_{10}$.}
  \label{fig:noncircuits:contracted}
\end{figure}

The edge set $c$ from Figure~\ref{fig:noncircuits:graph} is a circuit,
but each edge set $c_n$ it induces on a contracted graph $G_n$ with
$n=6k+4$ is not a circuit (Figure~\ref{fig:noncircuits:contracted}). Indeed, each $G_{6k+4}$ consists of
$G[V(T_{6k+4})]$, for each $i$ with $1\le i\le k$ a vertex
$v^{6k+4}_i$ corresponding to a contracted tail of the ladder $L_i$,
and a vertex $v^{6k+4}_0$ corresponding to the contracted tail of the
wide ladder and all ladders $L_j$ with $j>k$. The edge set $c_{6k+4}$
is not a circuit since it has degree $4$ at $v^{6k+4}_0$. Therefore,
$c$ is a circuit although it is the limit of the non-circuits
$c_{6k+4}$.

One can easily manipulate the example so that no $c_n$ with $n$ large
enough is a circuit by attaching copies $H_1,\dotsc,H_5$ of $G$ to $G$
by connecting the vertices of the first rung of the wide ladder in
$H_i$ to some suitable vertices of $L_i$.

\section{Singular homology}

A~more subtle approach than \v{C}ech homology, which has been pursued
in~\cite{Hom1}, is to see to what extent $\CCC(G)$ can be captured by
the singular homology of~$|G|$.  After all, $\CCC(G)$~was defined via
(the edge sets of) circles in~$|G|$, which are just injective singular
loops. Can we extend this correspondence between injective loops and
circuits to one between $H_1(|G|)$ (singular) and~$\CCC(G)$?

There are two things to notice about~$H_1(|G|)$. The first is that we
can subdivide a $1$-simplex (or concatenate two $1$-simplices into one
by the inverse procedure) by adding a boundary. Indeed, if
$\sigma\colon [0,1]\to |G|$ is a path in $|G|$ from $x$ to~$y$, say,
and $z$ is a point on that path, there are paths $\sigma'$ from $x$
to~$z$ and $\sigma''$ from $z$ to~$y$ such that
$\sigma'+\sigma''-\sigma$ is the boundary of a singular $2$-simplex
`squeezed' on to the image of~$\sigma$. The second fact to notice is
that inverse paths cancel in pairs: if $\sigma^+$ is an $x$--$y$ path
in~$|G|$, and $\sigma^-$ an $y$--$x$ path with the same image
as~$\sigma^+$, then $[\sigma^+ + \sigma^-] = 0\in H_1$.%
  \footnote{To see that this sum is a boundary, subtract the constant
    $1$-simplex $\sigma$ with value~$x$: there is an obvious singular
    $2$-simplex of which $\sigma^+ + \sigma^- - \sigma$ is the
    boundary. Subtracting~$\sigma$ is allowed, since $\sigma = \sigma
    + \sigma - \sigma$, too, is a boundary: of the constant
    $2$-simplex with value~$x$.}
These two facts together imply that every homology class in $H_1$ is
represented by a single loop: given any $1$-cycle, we first add pairs
of inverse paths between the endpoints of its simplices to make its
image connected in the right way, and then use Euler's theorem to
concatenate the $1$-simplices of the resulting chain into a single
loop~$\sigma$. Moreover, we may assume that this loop is based at a
vertex.

To establish the desired correspondence between $H_1(|G|)$
and~$\CCC(G)$, we would like to assign to a homology class
in~$H_1(|G|)$, represented by a single loop~$\sigma$, an edge set
$f([\sigma])\in\CCC(G)$. Intuitively, we do this by counting for each
edge $e$ of $G$ how often $\sigma$ traverses it entirely (which, since
the domain of $\sigma$ is compact, is a finite number of times), and
let $f([\sigma])$ be the set of those edges $e$ for which this number
is odd. Using the usual tools of homology theory, one can make this
precise in such a way that $f$ is clearly a well defined homomorphism
$H_1(|G|)\to\EEE(G)$,%
  \footnote{For each edge~$e$, let $f_e\colon |G|\to S^1$ be a map
    wrapping $e$ once round~$S^1$ and mapping all of $|G|\sm\kreis e$
    to one point of~$S^1$. Let $\pi$ denote the group isomorphism
    $H_1(S^1)\to\FF_2$. Given $h\in H_1(|G|)$, let $f(h)\assign
    \{\,e\mid (\pi\circ (f_e)_*)(h) = 1\in\FF_2\,\}$. See~\cite{Hom1}
    for details.}
and whose image is easily seen to be~$\CCC(G)$. What is not clear at
once is whether $f$ is $1$--$1$ and onto.

Surprisingly, $f$~is indeed surjective---and this is not even hard to
show. Indeed, let an edge set $D\in\CCC(G)$ be given. Our task is to
find a loop~$\sigma$ that traverses every edge in~$D$ an odd number of
times, and every other edge of $G$ an even number of times. As a first
approximation, we let $\sigma_0$ be a path that traverses every edge
of some fixed normal spanning tree of $G$ exactly twice, once in each
direction; see~\cite[Sec.~3.3]{RDsBanffSurvey} for how to construct
such a loop. Moreover, we construct $\sigma_0$ in such a way that it
pauses at every vertex~$v$---more precisely, so that
$\sigma_0^{-1}(v)$ is a union of finitely many closed intervals at
least one of which is non-trivial.  Next, we write $D$ as a thin sum
$D = \sum_i C_i$ of circuits; such a representation of $D$ exists by
definition of~$\CCC(G)$. For each of these $C_i$ we pick a vertex
$v_i\in\overline{C_i}$, noting that no vertex of $G$ gets picked
infinitely often, because it has only finitely many incident edges and
the $C_i$ form a thin family.  Finally, we turn $\sigma_0$ into the
desired loop~$\sigma$ by expanding the pause at each vertex $v$ to a
loop going once round every $\overline{C_i}$ with $v=v_i$. It is not
hard to show that $\sigma$ is continuous~\cite{Hom1}, and clearly it
traverses every edge of $G$ the desired number of times.\looseness=-1

Equally surprisingly, perhaps, $f$~is usually not injective (see
below). In summary, therefore, the topological cycle space $\CCC(G)$
of $G$ is related to the first singular homology group of~$G$ as
follows:

\begin{theorem}[\cite{Hom1}]\label{thm:Csingular}
  The map $f\colon H_1(|G|)\to\EEE(G)$ is a group homomorphism
  onto~$\CCC(G)$, which has a non-trivial kernel if and only if $G$
  contains infinitely many (finite) circuits.
\end{theorem}

An example of a non-null-homologous loop in~$|G|$ whose homology class
maps to the empty edge set $\es\in\CCC(G)$ is easy to describe. Let
$G$ be the one-way infinite ladder~$L$ (with its end on the right),
and define a loop~$\rho$ in~$L$, as follows. We start at time~$0$ at
the top-left vertex, $v_0$~say, and begin by going round the first
square of $G$ in a clockwise direction. This takes us back to~$v_0$.
We then move along the horizontal edge incident with~$v_0$, to its
right neighbour~$v_1$. From here, we go round the second square in a
clockwise direction, back to~$v_1$ and on to its right
neighbour~$v_2$. We repeat this move until we reach the end~$\omega$
of~$G$ on the right, say at time $\frac12\in [0,1]$. So far, we have
traversed the first vertical edge and every bottom horizontal edge
once (in the direction towards~$v_0$), every other vertical edge twice
(once in each direction), and every top horizontal edge twice in the
direction towards the end. From there, we now use the remaining half
of our time to go round the infinite circle formed by the first
vertical edge and all the horizontal edges one and a half times, in
such a way that we end at time~1 back at~$v_0$ and have traversed
every edge of~$L$ equally often in each direction. Clearly, $f$~maps
(the homology class of) this loop~$\rho$ to $0\in\CCC(L)$.

\begin{figure}[htbp]
  \centering
  \includegraphics{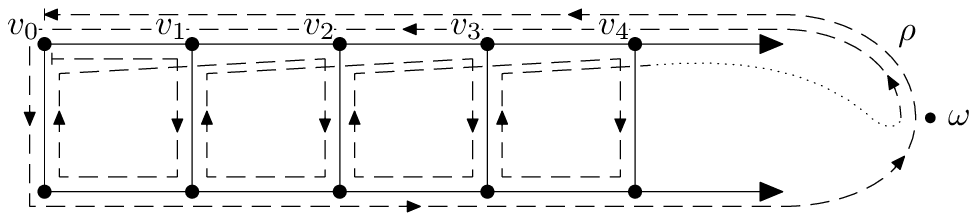}
  \caption{The loop $\rho$ is not null-homologous, but $f([\rho])
    =\es$.}
  \label{fig:kringel}
\end{figure}

The loop $\rho$ is indeed not null-homologous~\cite{Hom1}, but it
seems non-trivial to show this. To see why this is hard, let us
compare $\rho$ to a loop winding round a finite ladder in a similar
fashion, traversing every edge once in each direction. Such a
loop~$\sigma$ is still not null-homotopic, but it is null-homologous.
To see this, we subdivide it into single edges: we find a finite
collection of $1$-simplices~$\sigma_i$, four for every edge on the
topp and two for every other edge, such that $[\sigma] =
\big[\sum_i\sigma_i\big]$ and every $\sigma_i$ just traverses its
edge. Next, we pair up these~$\sigma_i$ into cancelling pairs: if
$\sigma_i$ and~$\sigma_j$ traverse the same edge~$e$ (in opposite
directions), then $[\sigma_i + \sigma_j] = 0$. Hence $[\sigma] =
\big[\sum_i\sigma_i\big] = 0$, as claimed. But we cannot imitate this
proof for $\rho$ and the infinite ladder, because homology classes in
$H_1(|G|)$ are still finite chains: we cannot add infinitely many
boundaries to subdivide $\rho$ infinitely often.

As it happened, the proof of the seemingly simple fact that $\rho$ is
not null-homologous took a detour via the solution of a much more
fundamental problem: the problem of understanding the fundamental
group of~$|L|$, or more generally, of~$|G|$ for a locally finite
graph~$G$. In order to distinguish $\rho$ from boundaries, we looked
for a numerical invariant $\Lambda$ of $1$-chains that was non-zero
on~$\rho$ but both linear and additive (so that
$\Lambda(\sigma_1\sigma_2) = \Lambda(\sigma_1 + \sigma_2) =
\Lambda(\sigma_1) + \Lambda(\sigma_2)$ for concatenations of
1-simplices~$\sigma_1,\sigma_2$) and invariant under homotopies (so
that $\Lambda(\sigma_1\sigma_2) = \Lambda(\sigma)$ when $\sigma\sim
\sigma_1\sigma_2$). Then, given a $2$-simplex $\tau$ with boundary
$\partial\tau = \sigma_1 + \sigma_2 - \sigma$, we would have
$\Lambda(\partial\tau) = \Lambda(\sigma_1\sigma_2) - \Lambda(\sigma) =
0$, so~$\Lambda$ would vanish on all boundaries but not on~$\rho$. We
did not quite find such an invariant~$\Lambda$, but a collection of
similar invariants which, together, can distinguish loops like~$\rho$
from boundaries.

\section{The fundamental group of $|G|$}

In this section we will sketch the combinatorial description
of~$\pi_1(|G|)$ given in~\cite{FundGp}. Our description involves
infinite words and their reductions in a `continuous' setting, and
embedding the group they form as a subgroup of a limit of finitely
generated free groups.

Let $G$ be a locally finite connected graph, fixed throughout this
section, and let $T$ be a topological spanning tree of $|G|$.
When $G$ is finite, then $\pi_1(|G|)=\pi_1(G)$~is the free group $F$
on the set $\{e_0,\dotsc,e_n\}$ of chords of any fixed spanning tree.
The standard description of $F$ is given in terms of reduced words of
those oriented chords, where reduction is performed by cancelling
adjacent inverse pairs of letters such as $\ve[i]\ev[i]$
or~$\ev[i]\ve[i]$. The map assigning to a path in $|G|$ the sequence
of (oriented) chords it traverses defines the canonical group
isomorphism between $\pi_1(|G|)$ and~$F$; in particular, reducing the
words obtained from homotopic paths yields the same reduced word.

Our description of $\pi_1(|G|)$ when $G$ is infinite is similar in
spirit, but more complex. We start not with an arbitrary spanning
tree but with a topological spanning tree of~$|G|$. Then every path
in $|G|$ defines as its `trace' an infinite word in the oriented
chords of that tree, as before. However, these words can have any
countable order type, and it is no longer clear how to define the
reduction of words in a way that captures homotopy of paths.

\begin{figure}[htbp]
  \centering
  \includegraphics[width=.7\textwidth]{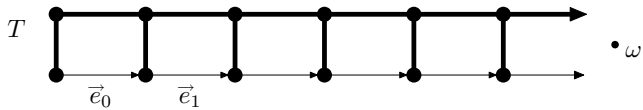}
  \caption{The infinite ladder and its topological spanning tree $T$
  (bold edges)}
  \label{fig:singleladder}
\end{figure}

Consider the following example. Let $G$ be the infinite ladder, with
a topological spanning tree~$T$ consisting of one side of the ladder,
all its rungs, and its unique end~$\omega$
(Figure~\ref{fig:singleladder}). The path running along the bottom
side of the ladder and back is a null-homotopic loop. Since it traces
the chords $\ve[0], \ve[1],\dotsc$ all the way to~$\omega$ and then
returns the same way, the infinite word $\ve[0]\ve[1]\dotso\ev[1]
\ev[0]$ should reduce to the empty word. But it contains no cancelling
pair of letters, such as $\ve[i]\ev[i]$ or~$\ev[i]\ve[i]$.

This simple example suggests that some transfinite equivalent of
cancelling pairs of letters, such as cancelling inverse pairs of
infinite sequences of letters, might lead to a suitable notion of
reduction. However, in graphs with infinitely many ends one can have
null-homotopic loops whose trace of chords contains no cancelling pair
of subsequences whatsoever: 

\begin{example}\label{ex:T2}
  We construct a locally finite graph $G$ and a null-homotopic loop
  $\sigma$ in~$|G|$ whose trace of chords contains no cancelling pair
  of subsequences, of any order type.

  Let $T$ be the binary tree with root $r$. Like in~\cite[pp.
  30--31]{RDsBanffSurvey} we can construct a loop $\sigma$ in $|T|$
  that traverses every edge of $T$ once in each direction, see
  Figure~\ref{fig:binary}.
  
  \begin{figure}[htbp]
    \centering
    \includegraphics{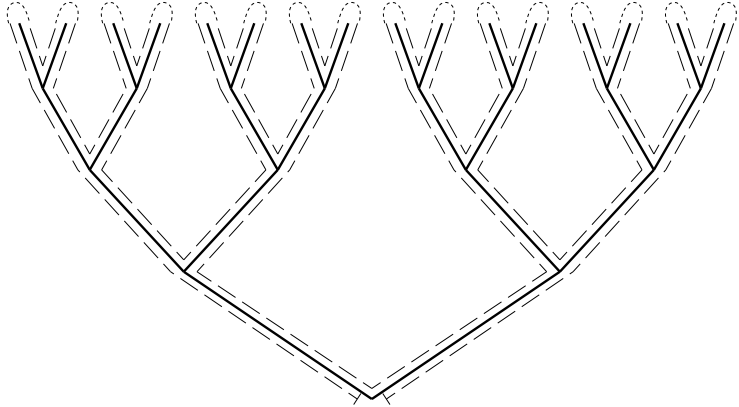}
    \caption{A loop running twice through each edge of the binary tree.}
    \label{fig:binary}
  \end{figure}
  
  The loop $\sigma$ is easily seen to be null-homotopic. It is also
  easy to check that no sequence of passes of $\sigma$ through the
  edges of $T$ is followed immediately by the inverse of this
  sequence.

  The edges of $T$ are not chords of a topological spanning tree, but
  this can be achieved by changing the graph: just double every edge.%
    \footnote{And subdivide the new edges once, in case you prefer to
    obtain a simple graph instead of a graph with multiple edges.}
  The new edges together with all vertices and ends then form a
  topological spanning tree in the resulting graph~$G$, whose chords
  are the original edges of our tree~$T$, and $\sigma$ is still a
  (null-homotopic) loop in~$|G|$.
\end{example}

Example~\ref{ex:T2} shows that there is no hope of capturing
homotopies of loops in terms of word reduction defined recursively by
cancelling pairs of inverse adjacent subwords, finite or infinite. We
shall therefore define the reduction of infinite words differently,
though only slightly. We shall still cancel inverse letters in pairs,
even only one at a time, and these reduction `steps' will be ordered
linearly (rather unlike the simultaneous dissolution of all the
chords by the homotopy in the example). However, the reduction steps
will not be well-ordered.

This definition of reduction is less straightforward, but it has an
important property: as for finite~$G$, it will be purely combinatorial
in terms of letters, their inverses, and their linear order, making no
reference to the interpretation of those letters as chords and their
relative positions under the topology of~$|G|$. 

Another problem, however, is more serious: since the reduction steps
are not well-ordered, it will be difficult to handle
reductions---e.g.\ to prove that every word reduces to a unique
reduced word, or that word reduction captures the homotopy of loops,
i.e.\ that traces of homotopic loops can always be reduced to the same
word. The key to solving these problems will lie in the observation
that the property of being reduced can be characterized in terms of
all the finite subwords of a given word. We shall formalize this
observation by way of an embedding of our group $F_\infty$ of
infinite words in the inverse limit $F^*$ of the free groups on the
finite subsets of letters.

\medbreak

A \emph{word} is a map
\begin{equation*}
  w \colon S \to A \assign \{\ve[0],\ve[1],\dotsc\} \cup
  \{\ev[0],\ev[1],\dotsc\}
\end{equation*}
(the letter $\ev[i]$ being the inverse of $\ve[i]$), where $S$ is a
totally ordered (countable) set, the set of \emph{positions} of (the
letters used by)~$w$, and every letter has only finitely many
preimages in $S$. A~\emph{reduction} of a word $w$ is a totally
ordered set $R$ of disjoint pairs of positions of~$w$ such that the
positions in each pair are mapped to inverse letters and are adjacent
in the word obtained from $w$ by deleting all (positions of) letters
contained in earlier pairs in $R$. We say that \emph{$w$ reduces to}
the word $w \restr (S\sm\bigcup R)$. If $w$ has no nonempty
reduction, we call it \emph{reduced}. Note that neither the set $S$
of positions of a word $w$ nor a reduction of $w$ have to be
well-ordered.

It was shown in~\cite{FundGp} that every word $w$ reduces to a unique
word $r(w)$%
  \footnote{Unique as an abstract word, not as a restiction of $w$:
  The word $\ve[0]\ev[0]\ve[0]$, for example, reduces to $\ve[0]$,
  but this letter can have the first or the last position in the
  original word}
and hence the reduced words form a group $F_{\infty}$. It was also
shown that $F_{\infty}$ embeds canonically in the inverse limit of
the groups $F_n$, the free groups on the sets $\{e_0,\dotsc,e_n\}$.

On the other hand, the fundamental group of $|G|$ embeds in
$F_{\infty}$: Mapping a homotopy class $\langle\alpha\rangle$ to the
word $r(w_{\alpha})$, where $w_{\alpha}$ is the \emph{trace} of
$\alpha$, the word induced by the passes of $\alpha$ through the
chords of $T$ (with their natural order given by $\alpha$), turns
out to be well-defined; in other words, the traces of homotopic loops
reduce to the same word. The harder part is to show the converse:
that two loops are homotopic whenever their traces reduce to the same
word. In~\cite{FundGp}, it was shown that the homotopy can even be
chosen so that it contracts pairs of passes, one at a time, like
known from finite graphs.

The map $\langle\alpha\rangle\mapsto r(w_\alpha)$ is not normally
surjective. For example, $\ve[0]\ve[1]\dotsm$ will always be a
reduced word, but no loop in $|G|$ can pass through these chords in
precisely this order if they do not converge to an end. Hence if
there is a non-converging sequence of chords---which is the case
whenever there are two ends of $G$ with no contractible neighbourhood
in $|G|$---then the reduced word $\ve[0]\ve[1]\dotsm$ lies outside
the image of our map $\langle\alpha\rangle\mapsto r(w_\alpha)$.

In order to describe the image of this map precisely, let us call a
word $w\colon S\to A$ \emph{monotonic} if there is an enumeration
$s_0,s_1,\dotsc$ of $S$ such that either $s_0<s_1<\dotsb$ or
$s_0>s_1>\dotsb$. Let us say that $w$ \emph{converges} if the
sequence of chords corresponding to its sequence $w_(s_0),w(s_1),
\dotsc$ of letters converges. If $w$ is the trace of a loop in~$|G|$,
then by the continuity of this path all the monotonic subwords of
$w$---and hence those of~$r(w)$---converge. It was shown
in~\cite{FundGp} that the converse is also true: A reduced word is
the trace of a loop in $|G|$ if and only if all its monotonic
subwords converge.

We can now summarize our combinatorial description of~$\pi_1(|G|)$ as
follows.

\begin{theorem}[\cite{FundGp}]\label{thm:pi1}
  Let $G$ be a locally finite connected graph, let $T$ be a
  topological spanning tree of~$|G|$, and let $e_0, e_1,\dots$ be its
  chords.
  \begin{enumerate}
  \item\label{enum:pi1Finf}
    The map $\langle\alpha\rangle\mapsto r(w_\alpha)$ is an injective
    homomorphism from $\pi_1(|G|)$ to the group $F_\infty$ of reduced
    finite or infinite words in $\{\ve[0],\ve[1],\dotsc\}\cup\{\ev[0],
    \ev[1],\dots\}$. Its image consists of those reduced words whose
    monotonic subwords all converge in~$|G|$.
  \item\label{enum:FinfProjLim}
    The homomorphisms $w\mapsto r(w\restr I)$ from $F_\infty$ to~$F_I$
    embed $F_\infty$ as a subgroup in~$\varprojlim F_I$. It consists
    of those elements of~$\varprojlim F_I$ whose projections $r(w
    \restr I)$ use each letter only boundedly often. (The bound may
    depend on the letter.)
  \end{enumerate}
\end{theorem}

\medskip

Theorem~\ref{thm:pi1} provides an interesting interaction between the
topological cycle space of $G$ and the fundamental group of $|G|$: It
is a well-known fact that the first (singular) homology group of a
space is the abelianization of its fundamental group. For graphs, this
yields that the (classical) cycle space of $G$ is the abelianization of
$\pi_1(G)$. Theorem~\ref{thm:pi1} implies an analoguous result for the
topological cycle space: It is the \emph{strong abelianization} of
$\pi_1(|G|)$~\cite[Theorem 6.19]{PhilippDiss}, the quotient of
$\pi_1(|G|)$ obtained by factoring out all words in which every letter
appears as often as its inverse.

\section{An ad-hoc homology for locally compact spaces}\label{sec:adhoc}

In this section we take up the thread of defining $\CCC(G)$ in terms
of homology. We have seen that \v{C}ech homology---although its first
group is isomorphic to the topological cycle space---fails to
properly reflect its relation to the combinatorial structure of $G$.
For this reason, we shall keep at our singular approach to define
$\CCC$ in terms of homology. Since by Theorem~\ref{thm:Csingular}
standard singular homology is not the right theory to capture $\vCCC$,
we shall define a singular-type homology that does so.

As advertised in Section~\ref{sec:intro}, we shall define our homology
for locally compact Hausdorff spaces with a (fixed) Hausdorff
compactification. Recall that these properties are needed to reflect
the properties of $G$ and $|G|$ that are fundamental for the success
of $\CCC$. Therefore, this class of spaces is the broadest for which
we can hope to obtain a homology theory with similar properties as
$\CCC$. Note that this class includes, for instance, all locally
finite CW-complexes, of any dimension.

Loops like the one in Figure~\ref{fig:kringel} suggest that our
homology should allow to subdivide a $1$-simplex infinitely often:
Then, every $1$-chain in $|G|$ will be homologous to the sum of its
passes through edges of $G$, and hence it will be null-homologous if
and only if it lies in the kernel of $f$. The idea is thus to define
the homology so that we obtain essentially the same $1$-cycles as in
standard singular homology but more boundaries.

The construction of $\CCC$ is based on the idea to consider not only
the graph itself but also its ends. Nevertheless, although ends do not
play a different role in the definition of $\CCC$ than points in $G$,
elements of $\CCC$ do behave differently at ends. Indeed, elements of
$\CCC$ are thin sums of circuits, and as $G$ is locally finite, these
circuits are also `thin' at vertices, i.e.\ every vertex lies in only
finitely many of the closures of the circuits in the family. This does
not have to be the case for ends: An end can lie in the closures of
infinitely many circuits, even when the circuits form a thin family.

This suggests to require a similar property from the chains in our
homology: They will have to be locally finite in $G$ but not at ends.%
  \footnote{The formal definition of `locally finite' will be given
  shortly.}
This will enable us to subdivide paths in $|G|$ infinitely often, but
the required locally finiteness in $G$ will keep us from obtaining
undesired cycles, such as the edges of a double-ray (all directed the
same way), which has zero boundary but does not correspond to an
element of the cycle space. In the ad-hoc homology we shall define in
this section we will rule out such cycles by imposing an additional
condition on cycles. This will lead to the desired result in
dimension $1$, i.e.\ our first homology group will be $\CCC$, but
generate problems elsewhere. More precisely, this homology will fail
to satisfy the Eilenberg-Steenrod axioms for homology, which is
caused precisely by this restriction on cycles.

In~\cite{Hom2} we thus change our approach slightly: Instead of
restricting the group of cycles we define chains differently, so as
to obtain $1$-cycles that are essentially finite and $2$-cycles that
allow us to subdivide $1$-simplices infinitely often. This homology
theory then satisfies the axioms~\cite{Hom2}. On the other hand, the
proof that this homology theory specializes in dimension $1$ to yield
$\CCC$ relies on the corresponding result for the ad-hoc homology
defined in this section. Moreover, it introduces some of the main
ideas from~\cite{Hom2} in a technically simpler setting.

\medskip

Let $X$ be a locally compact Hausdorff space and let $\hat X$ be a
Hausdorff compactification of $X$. (See e.g.~\cite{AbelsStrantzalos}
for more on such spaces.) Note that every locally compact Hausdorff
space is Tychonoff, and thus has a Hausdorff compactification.
Although we do not make any assumptions on the type of the
compactification, apart from being Hausdorff, we will call the points
in $\hat X\sm X$ \emph{ends}, even if they are not ends in the usual,
more restrictive, sense.

Let us call a family $(\sigma_i\mid i\in I)$ of singular $n$-simplices
in $\hat X$ \emph{admissible} if
\begin{enumerate}
\item\label{enum:locfininX}
  $(\sigma_i\mid i\in I)$ is locally finite in~$X$, that is, every
  $x\in X$ has a neighbourhood in $X$ that meets the image of
  $\sigma_i$ for only finitely many~$i$;
\item\label{enum:rootedinX}
  every $\sigma_i$ maps the $0$-faces of~$\Delta^n$ to~$X$.
\end{enumerate}
Note that as $X$ is locally compact, \ref{enum:locfininX} is
equivalent to asking that every compact subspace of $X$ meets the
image of $\sigma_i$ for only finitely many~$i$. Condition
\ref{enum:rootedinX}, like \ref{enum:locfininX}, underscores that ends
are not treated on a par with the points in~$X$: we allow them to
occur on infinitely many~$\sigma_i$ (which \ref{enum:locfininX}
forbids for points of $X$), but not in the fundamental role of images
of $0$-faces: all simplices must be `rooted' in~$X$.

When $(\sigma_i\mid i\in I)$ is an admissible family of $n$-simplices,
any formal linear combination $\sum_{i\in I} \lambda_i \sigma_i$ with
all $\lambda_i\in\ZZ$ is an \emph{$n$-sum in~$X$}.%
  \footnote{In standard singular homology, one does not usually
  distinguish between formal sums and chains. It will become apparent
  soon why we have to make this distinction.}
We regard $n$-sums $\sum_{i\in I}\lambda_i\sigma_i$ and
$\sum_{j\in J}\mu_j\tau_j$ as \emph{equivalent} if for every
$n$-simplex $\rho$ we have $\sum_{i\in I, \sigma_i=\rho}\lambda_i =
\sum_{j\in J, \tau_j=\rho}\mu_j$. Note that these sums are
well-defined since an $n$-simplex can occur only finitely many times
in an admissible family. We write $C_n(X)$ for the group of
\emph{$n$-chains}, the equivalence classes of $n$-sums. The elements
of an $n$-chain are its \emph{representations}. Clearly every
$n$-chain $c$ has a unique representation whose simplices are pairwise
distinct---which we call the \emph{reduced representation} of $c$---,
but we shall consider other representations too. The subgroup of
$C_n(X)$ consisting of those $n$-chains that have a finite
representation is denoted by $C'_n(X)$.

The boundary operators $\partial_n\colon C_n\to C_{n-1}$ are defined
by extending linearly from~$\partial_n\sigma_i$, which are defined as
usual in singular homology. Note that $\partial_n$ is well defined
(i.e., that it preserves the required local finiteness), and
$\partial_{n-1}\partial_n = 0$. Chains in $\im\partial$ will be called
\emph{boundaries}.

As $n$-cycles, we do \emph{not} take the entire kernel of
$\partial_n$. Rather, we define $Z'_n(X) \assign \ker(\partial_n\restr
C'_n(X))$, and let $Z_n(X)$ be the set of those $n$-chains that are
sums of such finite cycles:
\begin{equation*}
  Z_n (X) \assign \Big\{\varphi\in C_n(X)\Bigm| \varphi = \sum_{j\in J}
  z_j \text{ \em with } z_j\in Z'_n(X)\ \forall j\in J\Big\}.
\end{equation*}
More precisely, an $n$-chain $\varphi\in C_n(X)$ shall lie in $Z_n(X)$
if it has a representation $\sum_{i\in I}\lambda_i\sigma_i$ for which
$I$ admits a partition into finite sets~$I_j$ ($j\in J$) such that, for
every $j\in J$, the $n$-chain $z_j \in C'_n(X)$ represented by
$\sum_{i\in I_j} \lambda_i\sigma_i$ lies in $Z'_n(X)$. Any such
representation of $\varphi$ as a formal sum will be called a
\emph{standard representation} of~$\varphi$ \emph{as a cycle}.%
  \footnote{Since the $\sigma_i$ need not be distinct, $\varphi$~has many
  representations by formal sums. Not all of these need admit a
  partition as indicated---an example will be given later in the
  section.}
We call the elements of $Z_n(X)$ the \emph{$n$-cycles} of~$X$.

The chains in $B_n(X) \assign \im\partial_{n+1}$ then form a subgroup
of $Z_n(X)$: by definition, they can be written as
$\sum_{j\in J}\lambda_jz_j$ where each $z_j$ is the (finite) boundary
of a singular ($n+1$)-simplex. We therefore have homology groups
\begin{equation*}
  H_n(X) \assign Z_n(X)/B_n(X)
\end{equation*}
as usual.

Note that if $X$ is compact, then all admissible families and hence
all chains are finite, so the homology defined above coincides with
the usual singular homology. The characteristic feature of this
homology is that while infinite cycles are allowed, they are always of
`finite character': in any standard representation of an infinite
cycle, every finite subchain is contained in a larger finite subchain
that is already a cycle.

\medskip

\begin{figure}[htbp]
  \centering
  \includegraphics[width=0.7\linewidth]{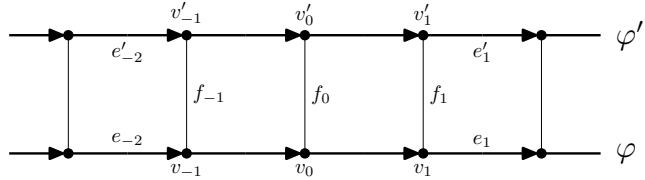}
  \caption{The $1$-chains $\varphi$ and $\varphi'$ in the double ladder.}
  \label{fig:double}
\end{figure}

Let us look at an example which might indicate whether we obtain the
desired cycles in order to capture the topological cycle space.
Consider the \emph{double ladder}. This is the $2$-ended graph $G$
with vertices $v_n$ and~$v'_n$ for all integers~$n$, and with edges
$e_n$ from $v_n$ to~$v_{n+1}$, edges $e'_n$ from $v'_n$ to~$v'_{n+1}$,
and edges $f_n$ from $v_n$ to $v'_n$. The $1$-simplices corresponding
to these edges, oriented in their natural directions, are
$\theta_{e_n}$, $\theta_{e'_n}$, and $\theta_{f_n}$, see Figure
\ref{fig:double}.

In order to let the elements of our homology be defined, let $\hat G$
be any Hausdorff compactification of $G$. (One could, for instance,
choose the Freudenthal compactification $|G|$ of $G$.) For the
infinite chains $\varphi$ and $\varphi'$ represented by $\sum
\theta_{e_n}$ and $\sum \theta_{e'_n}$, respectively, and for $\psi
\assign \varphi-\varphi'$ we have $\partial\varphi = \partial\varphi'
= \partial\psi = 0$, and neither sum as written above contains a
finite cycle. However, we can rewrite $\psi$ as $\psi = \sum z_n$ with
finite cycles $z_n = \theta_{e_n} + \theta_{f_{n+1}} - \theta_{e'_n} -
\theta_{f_n}$. This shows that $\psi\in Z_1(G)$, although this was not
visible from its original representation.

\begin{figure}[htbp]
  \centering
  \includegraphics[width=0.7\linewidth]{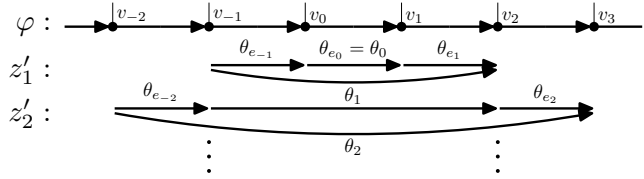}
  \caption{Finite cycles summing to~$\varphi$---by an inadmissible sum.}
  \label{fig:single}
\end{figure}

By contrast, one can show that $\varphi\notin Z_1(G)$ if $\hat G$ is
the Freudenthal compactification of $G$. This is proved
in~\cite{Hom1}, but is not obvious. For example, one might try to
represent $\varphi$ as $\varphi = \sum_{n=1}^{\infty}
z'_n$ with $z'_n \assign \theta_{e_{-n}} + \theta_{n-1} + \theta_{e_n}
- \theta_n$, where $\theta_n\colon[0,1]\to e_{-n}\cup\dotsb\cup e_n$
maps $0$ to $v_{-n}$ and $1$ to $v_{n+1}$, see
Figure~\ref{fig:single}.

This representation of $\varphi$, however, although well defined as a
formal sum (since every simplex occurs at most twice), is not a legal
$1$-sum, because its family of simplices is not locally finite and
hence not admissible. (The point $v_0$, for instance, lies in every
simplex $\theta_i$.)

\medskip

This homology indeed captures the cycle space~\cite{Hom1}. To see
this, note that since infinite chains are allowed, we can add
infinitely many boundaries to a loop like in Figure~\ref{fig:kringel}
so as to subdivide it into its edge passes. Note that the family of
boundaries we add has to be locally finite and it is not obvious that
this can always be satisfied. (See~\cite{Hom1} for how to choose the
boundaries.) Therefore, two chains are homologous if both of them
traverse each edge of $G$ the same number of times. Together with the
fact that the homomorphism $f$ from the first singular homology group
$H_1(|G|)$ to $\CCC(G)$ can be extended to a homomorphism $H_1(G) \to
\CCC(G)$~\cite{Hom1}, this implies that $H_1(G)$ and $\CCC(G)$ are
isomorphic.

\begin{theorem}[\cite{Hom1}]
  If $G$ is a locally finite graph and $\hat G=|G|$, then $H_1(G)$ is
  canonically isomorhic to $\CCC(G)$.
\end{theorem}

\medskip

Note that it does not suffice to require the chains to be locally
finite withour any further assumptions, as it is the case for the
\emph{locally finite homology} defined in~\cite{HughesRanicki}: This
homology does \emph{not} capture the cycle space. Indeed, applied to
$|G|$ it yields the usual singular homology, since every locally
finite chain in a compact space is finite. On the other hand, applied
to $G$, the locally finite homology allows for chains like $\varphi$
above, which do not correspond to an element of the cycle space.

\medskip

As mentioned before, the ad-hoc homology defined above does not
satisfy the Eilenberg-Steenrod axioms for homology. (For an example,
as well as a listing of the axioms, see~\cite{Hom2}.) This is caused
by the fact that the cycles are not chosen to be the entire kernel of
$\partial$ but with the additional property that they are a locally
finite sum of finite cycles.

For this reason, we develop in~\cite{Hom2} a homology that does
satisfy the axioms and that is defined without further assumptions
on the cycles. Like before, we define this homology for locally
compact Hausdorff spaces $X$ with a fixed Hausdorff compactification
$\hat X$. For this homology to capture $\CCC(G)$ we have to
allow infinite chains, since chains like (the chain consisting of)
the loop in Figure~\ref{fig:kringel} have to be null-homologous in
our homology---as they correspond to the empty edge set in
$G$---but are not the boundary of a finite chain. On the other hand,
we cannot allow all locally finite chains, as this would yield the
locally finite homology mentioned above. The solution to this
dilemma is surprisingly simple: We allow only those simplices to
appear infinitely often in a chain that are needed to subdivide a
path, or more generally, a simplex. This will enable us to
subdivide simplices into their edge passes and the isomorphism
between our new homology and $\CCC(G)$ will follow like for the
ad-hoc homology above.

A main feature of the simplices whose boundaries we need to
subdivide a path $\sigma$ is that they are in a sense
`one-dimensional': they can be written as the composition of a map
$\Delta^2\to\Delta^1$ and $\sigma$.%
  \footnote{Note that in general spaces the image of such a
  $2$-simplex does not have to be one-dimensional, since $\sigma$
  could be a space-filling curve.}
This leads us to the following definition: Call a singular
$n$-simplex $\tau$ in $\hat X$ \emph{degenerate} if there is a
compact Hausdorff space $X_{\tau}$ of topological dimension less
than $n$ such that $\tau$ can be written as the composition of
continuous maps $\Delta^n\to X_{\tau}\to \hat X$.

We would now like to say that we only allow chains (that have a
representation) with all but finitely many simplices degenerate.
This would not be a proper definition of `chain' since the boundary
of a chain would not have to be a chain in this case. This can
easily be remedied: Call a chain \emph{good} if it has the above
property. We now allow all $n$-chains that are the sum of a good
$n$-chain and the boundary of a good ($n+1$)-chain. This homology
turns out to satisfy all the Eilenberg-Steenrod axioms~\cite{Hom2},
and the fact that all $2$-simplices in the one-dimensional space
$|G|$ are degenerate implies that we indeed obtain the right
boundaries. Hence

\begin{theorem}[\cite{Hom2}]
  If $G$ is a locally finite graph and $\hat G=|G|$, then the first
  group $H_1(G)$ of the new homology is canonically isomorphic to
  $\CCC(G)$.
\end{theorem}

\bibliographystyle{amsplain}
\bibliography{collective}

\providecommand{\bysame}{\leavevmode\hbox to3em{\hrulefill}\thinspace}
\providecommand{\MR}{\relax\ifhmode\unskip\space\fi MR }
\providecommand{\MRhref}[2]{%
  \href{http://www.ams.org/mathscinet-getitem?mr=#1}{#2}
}
\providecommand{\href}[2]{#2}
\begin{thebibliography}{10}

\bibitem{AbelsStrantzalos}
H.~Abels and P.~Strantzalos, \emph{Proper transformation groups}, in
  preparation.

\bibitem{RDsBanffSurvey}
R.~Diestel, \emph{Locally finite graphs with ends: a topological approach}, \\
  {\tt http://arxiv.org/abs/0912.4213}, 2009.

\bibitem{TopSurveyI}
R.~Diestel, \emph{Locally finite graphs with ends: a topological approach.
  {I}.\ {B}asic theory}, Discrete Math. (to appear).

\bibitem{TopSurveyII}
R.~Diestel, \emph{Locally finite graphs with ends: a topological approach.
  {II}.\ {A}pplications}, Discrete Math. (to appear).

\bibitem{CyclesI}
R.~Diestel and D.~K{\"u}hn, \emph{On infinite cycles {I}}, Combinatorica
  \textbf{24} (2004), 68--89.

\bibitem{CyclesII}
\bysame, \emph{On infinite cycles {II}}, Combinatorica \textbf{24} (2004),
  91--116.

\bibitem{FundGp}
R.~Diestel and P.~Spr\"ussel, \emph{The fundamental group of a locally finite
  graph with ends}, arXiv \textbf{0910.5647} (2009).

\bibitem{Hom2}
\bysame, \emph{On the homology of locally compact spaces with ends}, arXiv
  \textbf{0910.5650} (2009).

\bibitem{Hom1}
\bysame, \emph{The homology of locally finite graphs with ends}, Combinatorica
  (to appear).

\bibitem{HughesRanicki}
B.~Hughes and A.~Ranicki, \emph{Ends of complexes}, Cambrigde Univ.\ Press,
  1996.

\bibitem{PhilippDiss}
P.~Spr\"ussel, \emph{On the homology of infinite graphs with ends}, Ph.D.
  thesis, Universit\"at Hamburg, 2010, see\\ {\small\tt
  http://www.sub.uni-hamburg.de/opus/volltexte/2010/4479/pdf/Dissertation.pdf}.

\end{thebibliography}

\small
\parindent=0pt
\vskip2mm plus 1fill

\begin{tabular}{cc}
\begin{minipage}[t]{0.5\linewidth}
Reinhard Diestel\\
Mathematisches Seminar\\
Universit\"at Hamburg\\
Bundesstra\ss e 55\\
20146 Hamburg\\
Germany\\
\end{minipage} 
&
\begin{minipage}[t]{0.5\linewidth}
Philipp Spr\"ussel\\
Mathematisches Seminar\\
Universit\"at Hamburg\\
Bundesstra\ss e 55\\
20146 Hamburg\\
Germany\\
\end{minipage}
\end{tabular}

\smallskip
Version 22.03.2010
\end{document}